\documentclass[final,12pt]{elsarticle}
\usepackage{xcolor,soul,cancel}
\usepackage[color]{showkeys}
\usepackage{amssymb}
\usepackage{amsmath}
\usepackage{amsthm}
\usepackage{graphicx}
\usepackage{float}
\journal{arXiv}

\definecolor{refkey}{rgb}{0,1,1}
\definecolor{labelkey}{rgb}{1,0,0}

\newtheorem{thm}{Theorem}

\newtheorem{prop}{Proposition}

\newcommand{\abs}[1]{\left\vert#1\right\vert}
\newcommand{\rk}{\operatorname{rank}}
\newcommand{\trace}{\operatorname{tr}}
\newcommand{\Span}{\operatorname{Span}}
\newcommand{\re}{\operatorname{Re}} 
\newcommand{\im}{\operatorname{Im}}

\newcommand{\eq} [1] {\begin{equation}\label{#1}\quad}
\newcommand{\en} {\end{equation}}
\newcommand{\scal}[1]{\langle#1\rangle}
\newcommand{\norm}[1]{\left\Vert#1\right\Vert}
\newcommand{\C}{\mathbb C}
\newcommand{\M}{{\bf M}}
\newcommand{\diag}{\operatorname{diag}}

\begin{document}

\begin{frontmatter}

\title{On low-dimensional partial isometries
\tnoteref{support}}

\author[nyuad]{Qixiao He}
\ead{qh433@nyu.edu}

\author[nyuad]{Ilya M. Spitkovsky}
\ead{ims2@nyu.edu, ilya@math.wm.edu, imspitkovsky@gmail.com}
\address[nyuad]{Division of Science and Mathematics,
New York  University Abu Dhabi (NYUAD), Saadiyat Island,
P.O. Box 129188 Abu Dhabi, United Arab Emirates}
\author[nyuad]{Ibrahim~Suleiman}
\ead{is1647@nyu.edu}
\tnotetext[support]{The results are partially based on the Capstone
projects of [QH] (2019-20) and [IS] (2021-22 academic year) under the supervision of [IMS]. The latter was also
supported in part by Faculty Research funding from the Division of Science
and Mathematics, New York University Abu Dhabi. }
\begin{abstract} Two statements concerning $n$-by-$n$ partial isometries are being considered: (i) these matrices are generic, if unitarily irreducible, and (ii) if nilpotent, their numerical ranges are circular disks. Both statements hold for $n\leq 4$ but fail starting with $n=5$. 
\end{abstract}

\end{frontmatter}

\section{Introduction}
Let $\M_n$ stand for the algebra of all $n$-by-$n$ matrices with the entries in the field $\C$ of complex numbers. The {\em numerical range} $W(A)$ of $A\in\M_n$
is defined as
\eq{nr} W(A)=\{ \scal{Ax,x}\colon x\in\C^n, \norm{x}=1 \}, \en
where $\norm{.}$ is the norm associated with the standard scalar product $\scal{.,.}$ on the space $\C^n$. It is well known that $W(A)$ is a compact (easy), convex (classical Toeplitz-Hausdorff theorem) subset of $\C$, invariant under unitary similarities of $A$. The shape of $W(A)$ is determined completely by the {\em Kippenhahn polynomial} of $A$, which by one of the currently accepted definitions (see, e.g., \cite[Chapter 6]{GauWu}) is the homogeneous polynomial $\det(x\re A+y\im A+zI)$. 
Here and in what follows, the standard notation $\re X=(X+X^*)/2$ and $\im X =i(X^*-X)/2$ is used for any $X\in\M_n$.

For computational purposes we prefer the {\em dehomogenized version}, which is the characteristic polynomial $P_{A,\theta}$ of $\re(e^{-i\theta}A)$. The roots of the latter, i.e. the eigenvalues 
\eq{lam} \lambda_1(\theta)\geq\lambda_2(\theta)\ldots\geq\lambda_n(\theta) \en of $\re(e^{-i\theta}A)$, define a family of lines 
\[ e^{i\theta}(\lambda_j(\theta)+i\mathbb R), \quad j=1,\ldots,n; \ \theta\in (-\pi,\pi], \]
the envelope of which is an algebraic curve $C(A)$ called the {\em Kippenhahn curve} of $A$. As it happens, $W(A)$ is the convex hull of $C(A)$. 

For our purposes, the following two facts are important:
	
{\sl Fact 1.} $C(A)$ contains a circle of radius $r$ centered at the origin if and only if $P_{A,\theta}$ is divisible by $\lambda^2-r^2$, and

{\sl Fact 2.} The boundary $\partial W(A)$ of the numerical range of $A$ itself is an analytic algebraic curve (and thus does not contain any line segments) if $A$ is generic. 

The term {\em generic} means, by definition (see, e.g., \cite[Definition 9]{JAG98}), that the inequalities in \eqref{lam} are strict for all $\theta\in (-\pi,\pi]$.

We refer the reader to a recent comprehensive monograph \cite{GauWu} for these and other properties of the numerical range.

Generic matrices are not normal; for $n=2$ the converse is also true. Namely, for $A\in\M_2$ the numerical range is an elliptical disk with the foci at the eigenvalues $\lambda_1,\lambda_2$ of $A$,  degenerating into the line segment $[\lambda_1,\lambda_2]$ if and only if $A$ is unitarily reducible (which in $n=2$ case is equivalent to $A$ being normal). Consequently, when $n=2$ it is true that

\begin{itemize}
\item[(NF)] Unitarily irreducible $A$ is generic (and thus the boundary $\partial W(A)$ of its numerical range does not contain any line segments), \end{itemize} and
\begin{itemize} \item[(Circ)] If $A$ is nilpotent, then $W(A)$ is a circular disk.
\end{itemize}

Both (NF) and (Circ) fail already for $n=3$. In particular, there exist unitarily irreducible nilpotent matrices $A\in\M_3$ with $\partial W(A)$ containing a flat portion; see \cite{KRS} for the complete description of such matrices.

In this paper, we determine the extent to which statements (NF) and (Circ) hold in the special case of {\em partial isometries}. Recall therefore that $A$ is a partial isometry if it preserves norms of the vectors in the orthogonal complement of its kernel:
\eq{pi} \norm{Ax}=\norm{x} \text{ for all } x\perp\ker A .\en
Equivalently, $A^*$ is an inner inverse of $A$, i.e. $AA^*A=A$, in which case it is actually the Moore-Penrose inverse $A^\dagger$ of $A$ (see, e.g., the survey \cite{GaPaRo} for these and other known facts about partial isometries).

We will always require $\ker A\neq\{0\}$ because otherwise $A$ is unitary and, as such, both unitarily reducible and invertible.

Note that, whenever (Circ) holds, the disk $W(A)$ is centered at the origin, because for any $A\in\M_n$ the foci of $C(A)$ coincide with the spectrum $\sigma(A)$ of $A$ while in our case $\sigma(A)=\{0\}$. 

If a partial isometry $A$ is not nilpotent but $W(A)$ is nevertheless a circular disk, it was conjectured in \cite{GWW} (and proved for $n\leq 4$) that it still has to be centered at the origin. This conjecture was proved for $n=5$ in \cite{SuSWe}, and for any $n$ under the additional requirement $\rk A=n-1$ in \cite{WeSp}, but remains unsolved in general. 

As in \cite{GWW,SuSWe}, we will make use of the block matrix representation of partial isometries as
\eq{BC} A=\begin{bmatrix}0 & B \\ 0 & C\end{bmatrix},  \en
where $B^*B+C^*C=I$. Applying an appropriate block diagonal unitary similarity, it is possible to put $C$ in a triangular form or to replace $B$ with the middle factor from its singular value decomposition. We will take advantage of this ability when convenient.

The paper is organized as follows. In Section~\ref{s:r} it is shown that both (NF) and (Circ) hold for partial isometries of rank one and $n-1$, implying in particular their validity for $n=3$. Section~\ref{s:r2} is devoted to 4-by-4 partial isometries. The validity of (NF) and (Circ) is established by way of considering the remaining case of rank two matrices. The dimension is again increased by one in Section~\ref{s:r3} where it is shown that (Circ) finally fails. More specifically, a criterion is established for 5-by-5 nilpotent partial isometries of rank three to have a circular numerical range, thus completely describing the (non-empty) set of those which do not. Within the latter, in Section~\ref{s:r3nf} we pinpoint a much smaller subset of matrices for which (NF) also fails. A short Section~\ref{s:hr} is about higher rank numerical ranges of the matrices considered in Sections~\ref{s:r2} and \ref{s:r3}.    

\section{Extreme rank values}\label{s:r}
Consider first the case of $A\in\M_n$ with $n$ arbitrary but $\rk A=1$. Then $A$ is unitarily similar to the direct sum of some $A_0\in\M_2$ with an $(n-2)$-dimensional zero block. From here it immediately follows

\begin{prop}\label{th:r1} If $A$ is a non-normal matrix of rank one, then $W(A)$ is an elliptical disk with the foci $0$ and
$\trace A$. This disk is circular if and only if $A$ is nilpotent. \end{prop}
In particular, (NF) and (Circ) hold in this case, whether or not $A$ is a partial isometry.

As it happens, the other extreme is also easy to handle.

\begin{prop}\label{th:n-1}
Let $A\in\M_n$ be a partial isometry with $\rk A=n-1$. Then both {\em (NF)} and {\em (Circ)} hold.
\end{prop}
\begin{proof} According to \cite[Proposition 2.3]{GWW} such $A$, if in addition it is unitarily irreducible, belongs to the so called class $S_n$. In other words, $A$ is a contraction with all its eigenvalues lying in the unit disk, and $I-A^*A$ has rank one. The hermitian parts $\re A$ of such matrices have only simple eigenvalues \cite[Corollary 2.7]{GauWu98}. Since the $S_n$ class is invariant under multiplication by unimodular scalars, this implies that all the eigenvalues of $\re(e^{-i\theta}A)$ are simple for any $\theta$. In other words, $S_n$ consists of generic matrices. In particular, every point of $\partial W(A)$ is an extreme point of $W(A)$ by \cite[Theorem 2.2]{GauWu13}. This takes care of (NF).

(Circ) follows immediately from the observation that a partial isometry with one-dimensional kernel is nilpotent if and only if it is unitarily similar to a Jordan block. \end{proof}

For $n=3$ Propositions~\ref{th:r1} and \ref{th:n-1} cover all the possibilities. Therefore, the following result holds.
\begin{thm}\label{th:n3}Statements {\em (NF)} and {\em (Circ)} hold for $3$-by-$3$ partial isometries.  \end{thm}
Observe that a unitarily irreducible rank 2 partial isometry $A\in\M_3$ is unitarily similar to
\eq{n3r2} \begin{bmatrix}0 & \sqrt{1-\abs{\lambda_1}^2} & -\overline{\lambda_1} \sqrt{1-\abs{\lambda_2}^2} \\
0 & \lambda_1 & \sqrt{1-\abs{\lambda_1}^2} \sqrt{1-\abs{\lambda_2}^2} \\ 0 & 0 & \lambda_2\end{bmatrix}, \en
where $\abs{\lambda_1},\abs{\lambda_2}<1$. This is a particular (for $n=3$) case of \cite[Theorem 1.4 in Chapter 7]{GauWu}, as well as the result of using the Takenaka-Malmquist basis (as described in \cite[Chapter 9]{DGSV}), but can also be verified directly.

Tests for possible shapes of the numerical ranges for $3$-by-$3$ matrices are well known; see \cite[Section 6.2]{GauWu} and references therein. When applied to \eqref{n3r2}, they yield the following description.
\begin{prop} \label{th:elvo}
The numerical range of the matrix \eqref{n3r2} is an elliptical disk with the foci $\lambda_1,\lambda_2$ if $\lambda_1=-\lambda_2$, 0 and $\lambda_1$ if $\lambda_1=\lambda_2$, and has an ovular shape otherwise.
\end{prop}
\begin{proof} 
Since the matrix \eqref{n3r2} is triangular, for the elliptical shape it is convenient to use the criterion as stated in \cite{KRS}, Theorems~2.2 and 2.4. According to them, for a unitarily irreducible $3$-by-$3$ matrix 
\eq{mat3by3} A= \begin{bmatrix}
  a & x & y\\
  0 & b & z\\
  0 & 0 & c
\end{bmatrix},\en
its numerical range is an elliptical disk if and only if \eq{lambda} \lambda := (c|x|^2 + b|y|^2 + a|z|^2 - x\overline{y}z)/(|x|^2 + |y|^2 + |z|^2) ,\en coincides with one of the eigenvalues $a,b,$ or $c$. The foci of the ellipse in question are then the other two eigenvalues. 

A direct computation shows that for the matrix \eqref{n3r2} as $A$, \eqref{lambda} takes the form 
\eq{lambda1} \lambda=\frac{\lambda_1(1-\abs{\lambda_2}^2)+\lambda_2(1-\abs{\lambda_1}^2)}{(1-\abs{\lambda_1}^2)+(1-\abs{\lambda_2}^2)}.\en
So, $\lambda$ is a convex combination of $\lambda_1$ and $\lambda_2$ with strictly positive coefficients. It therefore coincides with either of $\lambda_{1}$ or $\lambda_2$ if and only if $\lambda_1$ and $\lambda_2$ coincide. The foci of $W(A)$ are then $\lambda$ and the remaining eigenvalue, which is zero. 

On the other hand, the right hand side of \eqref{lambda1} is zero if and only if $\lambda_1=-\lambda_2$. If this is the case, the foci of $W(A)$ are $\lambda_1$ and $\lambda_2$. This completes the proof of the ellipticity criterion. 

The shape of $W(A)$ with a flat portion on the boundary is not an option, due to Proposition~\ref{th:n-1}. For completeness, however, we will provide  a dimension-specific independent proof, based on \cite[Theorem 2.11 in Chapter 6]{GauWu} (going back to \cite[Theorem~1.2]{RS05}). 

Namely, for a flat portion on the boundary of $W(A)$ of the matrix \eqref{mat3by3} to exist, it is necessary that $xyz\neq 0$ and 
\[ \left|\frac{xy}{z}\right| - 2\re(e^{-i\theta}a) = \left|\frac{xz}{y}\right| - 2\re(e^{-i\theta}b), \]
where $\theta = \arg(x\overline{y}z)$.
For $A$ as in \eqref{n3r2}, this condition can be rewritten as \[ \abs{\lambda_1}=\frac{1-\abs{\lambda_1}^2}{\abs{\lambda_1}}+2\abs{\lambda_1}, \]
which is obviously false. So indeed, when the numerical range is not an elliptical disk, i.e., $\lambda_1\neq\pm\lambda_2$, it is forced to be of an oular shape.  \end{proof}

\begin{figure}[H]
\centering
    \includegraphics[scale=0.5]{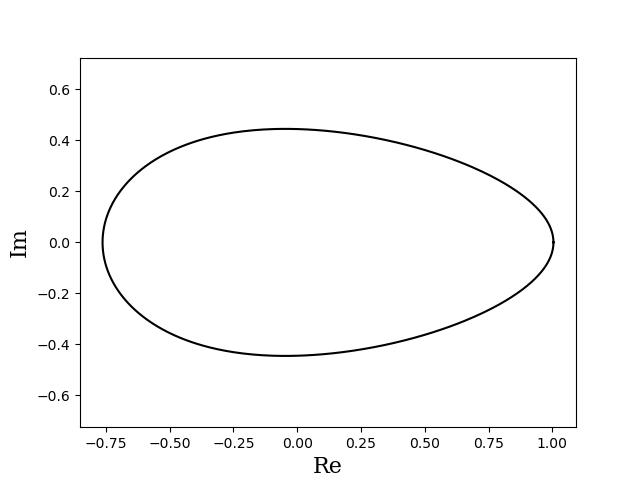}
    \caption{$W(A)$ for $A$ as in \eqref{3by3-ovular}}
    \label{ovularfig}
\end{figure}

We exhibit below a 3-by-3 partial isometry $A$ of rank two for which $W(A)$ is ovular. Let
\eq{3by3-ovular} A = \begin{bmatrix}
0 & \frac{1}{2} \, \sqrt{\frac{1}{2}} & -\frac{1}{8} \, \sqrt{15} \sqrt{\frac{7}{2}} \\
0 & \frac{1}{2} \, \sqrt{\frac{7}{2}} & \frac{1}{8} \, \sqrt{15} \sqrt{\frac{1}{2}} \\
0 & 0 & \frac{1}{4}
\end{bmatrix}.\en
Figure \ref{ovularfig} shows the ovular shape of $W(A)$.

\section{4-by-4 matrices} \label{s:r2}
Due to Propositions~\ref{th:r1} and \ref{th:n-1}, for $A\in\M_4$ we need only to consider the case $\rk A=2$.
\begin{prop}\label{th:nfn4r2}Unitarily irreducible $4$-by-$4$ partial isometries of rank two are generic. \end{prop}
\begin{proof} Let $A\in\M_4$ be a rank two partial isometry. Using the form of \eqref{BC} with the block $B$ replaced by the middle factor of its singular value decomposition, we may represent $A$ as  
\eq{An4r2} \begin{bmatrix} 0 & 0 & b_1 & 0 \\ 0 & 0 & 0 & b_2 \\ 0 & 0 & c_{11} & c_{12} \\  0 & 0 & c_{21} & c_{22} \end{bmatrix}.\en
Here $b_1,b_2\geq 0$, the columns $[c_{1j},c_{2j}]^T$ have lenghts $\sqrt{1-b_j^2}$ and are orthogonal to each other, $j=1,2$.

Applying a transposition $(1342)$ to the rows and columns of \eqref{An4r2}, observe that it is unitarily similar to the tridiagonal matrix
\eq{tri} \begin{bmatrix} 0 & b_1 & 0 & 0 \\ 0  & c_{11} & c_{12} & 0 \\ 0 & c_{21} & c_{22} & 0  \\ 0 & 0& b_2 & 0   \end{bmatrix}. \en
Suppose that the matrix \eqref{An4r2} (and thus also \eqref{tri}) is not generic. From \cite[Theorem 10]{BS041} (see also \cite[Lemma 4.12 in Chapter 6]{GauWu}) it then follows that at least one of the off-diagonal pairs in \eqref{tri} has coinciding absolute values. In other words, $b_1=0$, $b_2=0$, or $\abs{c_{12}}=\abs{c_{21}}$.

In the first two cases the unitary reducibility of $A$ is obvious. So, let us concentrate on the latter case. Along with the orthogonality of the columns of $C$ it implies that $\abs{c_{11}}=\abs{c_{22}}$, which makes $C$ a scalar multiple of the unitary matrix, and thus unitarily reducible. The unitary reducibility of $A$ then follows from \cite[Proposition 2.6]{GWW}. 

Note that it can also be checked directly that the span of $[0,\xi,\eta,0]^T$ and $[\xi,0,0,\eta]^T$, where $[\xi,\eta]^T$ is an eigenvector of $C$, forms a 2-dimensional, and thus non-trivial, reducing subspace of the matrix \eqref{tri}.

Since $A$ is given to be unitarily irreducible, the contradiction obtained completes the proof. \end{proof}
We move now to the nilpotent setting.
\begin{prop}\label{th:nilp42}Let $A\in\M_4$ be a nilpotent partial isometry of rank two. Then $W(A)$ is a circular disk.  \end{prop}
\begin{proof}Let us use representation \eqref{BC}, this time with $C$ put in an upper triangular form. In addition, let us use the block diagonal unitary similarity $\diag[U,I]$ to set the left lower entry of $B$ to zero. Then $A$ becomes
\eq{nil42} \begin{bmatrix} 0 & 0 & 1 & 0 \\ 0  & 0 & 0 & b \\ 0 & 0 & 0 & c  \\ 0 & 0& 0 & 0\end{bmatrix},\en
where without loss of generality $b,c\geq 0$, and $b^2+c^2=1$. Thus
\eq{nil42poly} P_{A,\theta}(\lambda) = \left(\lambda^2 - \frac{1}{2} + \sqrt{\frac{1}{4} - \frac{b^2}{16}}\right)\left(\lambda^2 - \frac{1}{2} - \sqrt{\frac{1}{4} - \frac{b^2}{16}}\right).\en 
According to Fact 1, the Kippenhahn curve $C(A)$ consists of the two circles centered at the origin of the radii \eq{r12} r_{1,2}=\sqrt{\frac{1}{2}\pm\sqrt{\frac{1}{4} - \frac{b^2}{16}}}.\en Consequently, $W(A)$ is a circular disk of radius $r_1$. \end{proof}
Taking into consideration Propositions~\ref{th:nfn4r2} and \ref{th:nilp42}, we arrive at the following
\begin{thm}\label{th:n4}Statements {\em (NF)} and {\em (Circ)} hold for $4$-by-$4$ partial isometries.  \end{thm}

\section{The statement (Circ) for 5-by-5 matrices} \label{s:r3} 
Due to Propositions~\ref{th:r1} and \ref{th:n-1}, we only need to consider matrices of rank two and three. Furthermore, rank two matrices $A\in\M_5$ are unitarily reducible, and so for them the (NF) statement is vacuously correct. (Circ) also holds, but for a different reason: a nilpotent partial isometry $A\in\M_5$ of rank two is unitarily similar to the direct sum of two nilpotent partial isometries of smaller size. The numerical ranges of the latter are concentric circular disks, and so $W(A)$ is just the larger of them.

So, we should concentrate on partial isometries $A\in\M_5$ of rank three. 

Via an appropriate unitary similarity and a rotation, a nilpotent rank three partial isometry $A\in\M_5$  can be put in the form
\eq{n5r3} \begin{bmatrix} 0 & 0 & 1 & 0 & 0 \\ 0 & 0 & 0 & b & tc \\
0 & 0 & 0 & c & -tb \\ 0 & 0 & 0 & 0 & s \\ 0 & 0 & 0 & 0 & 0 \end{bmatrix}, \en  
where \eq{bcdt} b,c,s,t \geq 0 \text{ and }s^2+t^2= b^2+c^2=1.\en  
	
\begin{prop}\label{th:cir5}The numerical range of the matrix \eqref{n5r3} is a circular disk if and only if \eq{bcdx} bcst=0.\en  If this is the case, then the radius of this disk (i.e., the numerical radius of $A$) is 
\[ r  = \frac{1}{2}\sqrt{\frac{3 + \sqrt{5-4(b^2+c^2t^2)}}{2}}. \] 
Moreover, $A$ is unitarily irreducible if and only if \eqref{bcdx} fails, or if exactly one of $b$ and $t$ is equal to zero.
\end{prop}	
It is not surprising that the numerical range of a unitarily reducible matrix \eqref{n5r3} is a circular disk: being nilpotent, it follows from the validity of (Circ) for $n\leq 4$. Observe though that in the unitarily irreducible case both circular and non-circular numerical ranges materialize.
\begin{proof}A direct computation shows that the Kippenhahn polynomial of the matrix \eqref{n5r3} is 	
\eq{n5r3kip-poly} P_{A,\theta}(\lambda) = -\lambda^{5}  + \frac{3}{4}\lambda^{3} - \frac{1}{16} \left(c^{2} t^{2} + b^{2} + 1\right)\lambda - \frac{1}{16}  b c{s} t \re (e^{i\theta}) ,\en
the homogeneous form of which being
\[ z^5-\frac{3}{4}(x^2+y^2)z^3+\frac{1}{16}(c^2t^2+b^2+1)(x^2+y^2)^2z-\frac{1}{16}bc{s}tx (x^2+y^2)^2. \]

If \eqref{bcdx} fails, the latter polynomial is irreducible, the Kippenhahn curve $C(A)$ does not contain circular components, and so $W(A)$ cannot possibly be a circular disk. Also, the irreducibility of the Kippenhahn polynomial implies unitary irreducibility of the matrix.

On the other hand, under condition \eqref{bcdx} \[ P_{A,\theta}=-\lambda\left(\lambda^2 - \frac{3}{8} + \frac{1}{8}\sqrt{5-4(b^2+c^2t^2)}\right)\left(\lambda^2 - \frac{3}{8} - \frac{1}{8}\sqrt{5-4(b^2+c^2t^2)}\right)\]
and so $C(A)$ is the union of the origin and two circles centered at the origin, having radii
\eq{rpm} r_\pm = \frac{1}{2}\sqrt{\frac{3 \pm \sqrt{5-4(b^2+c^2t^2)}}{2}}.\en

Consequently, $W(A)$ is indeed a circular disk of radius $r_+$. 

It remains to treat the unitary (ir)reducibility of $A$ provided that \eqref{bcdx} holds. 

To this end, observe that due to \eqref{bcdt} at most one of the variables $b,c$ can equal zero; the same is true for the pair $s,t$. Therefore, there are two cases to consider: 

{\sl Case 1.} Exactly two of the variables $b,c,s,t$ are equal to zero.

\noindent 
The possibilities are as follows: (i) $b=t=0, c=s=1$, (ii) $c=t=0, b=s=1$, (iii) $b=s=0, c=t=1$, and (iv) $c=s=0, b=t=1$. In each of the possibilities (i)--(iv), a simple permutational similiarity shows that $A$ is unitarily reducible.  

Before moving to the remaining situation, observe that the kernel of the skew-hermitian part of $A$ (which is a priori non-trivial, $A$ being a real matrix of odd size) is one-dimensional: $\ker\im(A) =\Span\{\xi\}$, where $\xi=[t,-s,0,tc,-b]^T$.

If $\mathcal L$ is a reducing subspace of $A$, then so is its orthogonal complement $\mathcal L^\perp$, and $\xi$ has to lie in one of them.  Switching the notation if needed, let $\xi\in\mathcal L$. Then also $A^k\xi, A^{*k}\xi\in\mathcal L$, $k=1,2,\ldots$

{\sl Case 2.} Exactly one of the variables $b,c,s,t$ is equal to zero. Again, there are four possibilities: (i) $b=0, c=1$, $s,t\neq 0$,\ \ (ii) $b=1, c=0$, $s,t\neq 0$, \ \ (iii) $t=0, s=1$, $b,c\neq 0$,  and (iv) $t=1, s=0$, $b,c\neq 0$. 

A straightforward computation shows that $\xi,A\xi,A^2\xi,A^{*2}\xi,A^{*3}\xi$ are linearly independent in case (i), and 
$\xi,A\xi,A^2\xi,A^{3}\xi,A^{*2}\xi$ are linearly independent in case (iii). So, in these cases $\mathcal L$ is the whole space $\C^5$ proving that $A$ is unitarily irreducible. 

On the other hand, in cases (ii),(iv) $A^3=0$, and so the non-trivial subspace $\mathcal L =\Span\{\xi,A\xi,A^2\xi\}$ is invariant under $A$. Another direct computation shows that $A^*A\xi,A^*A^2\xi\in\mathcal L$ (while $A^*\xi= A\xi$ by the choice of $\xi$), and so $\mathcal L$ is also invariant under $A^*$. This proves unitary reducibility of $A$ in cases (ii) and (iv) thus completing the proof. \end{proof} 

Note that Proposition~\ref{th:nilp42} and the criterion \eqref{bcdx} from Proposition~\ref{th:cir5} can also be established by using results from \cite{Mat10} (Remark~4 and Theorem~2, respectively).

\begin{figure}[H]
\centering
    \includegraphics[scale=0.5]{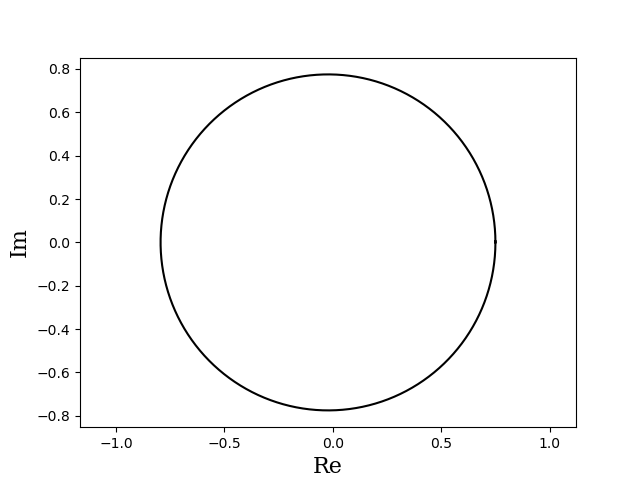}
    \caption{The plot of $W(A)$ looks circular, but in fact it is not.}
    \label{fig:noncircular}
\end{figure}

So, (Circ) fails even for partially isometric matrices, starting with $n=5$. A concrete example is given by
 \[ A = \begin{bmatrix} 0 & 0 & 1 & 0 & 0 \\ 0 & 0 & 0 & 1/2 & \ \sqrt{3}/{4} \\ 0 & 0 & 0 & {\sqrt{3}}/{2} & -{1}/{4} \\ 0 & 0 & 0 & 0 &  {\sqrt{3}}/{2} \\ 0 & 0 & 0 & 0 & 0\end{bmatrix}.\] See Figure \ref{fig:noncircular} for a plot of $W(A)$. Although the plot looks circular, the extreme eigenvalues of $\re(A)$ are $-0.79435\dots$ and 0.75 while those of $\im(A)$ are $\pm 0.77482\dots$. This shows that $W(A)$ is not circular.  

\section{The statement (NF) for 5-by-5 matrices} \label{s:r3nf} 

As it happens, (NF) fails when $n=5$ as well, and respective examples can also be found among nilpotent matrices. In this setting, up to unitary similarity and rotation, there are exactly two such matrices, as the next theorem shows. 
\begin{thm}\label{th:gen5}A nilpotent rank-three partial isometry $A\in\M_5$ is generic unless it is unitarily similar to 
\eq{spec}  e^{i\phi}\begin{bmatrix} 0 & 0 & 1 & 0 & 0 \\ 0 & 0 & 0 & \sqrt{1-c^2}  & \frac{c}{\sqrt{2\mp c}} \\
	0 & 0 & 0 & c & -\sqrt{\frac{1-c^2}{2\mp c}} \\ 0 & 0 & 0 & 0 & \sqrt{\frac{1\mp c}{2\mp c}}  
	\\ 0 & 0 & 0 & 0 & 0 \end{bmatrix}, \en 
where $\phi\in [0,2\pi]$ and the upper/lower choice of the sign corresponds to $c=c_\pm$ with  
\eq{cp} c_+:= \frac{2}{3}\left(1 - \sqrt{7}\cos\left(\alpha +\frac{\pi}{3} \right)\right)  = 0.55495\dots, \en 
\eq{cm} c_-:= -\frac{2}{3}\left(1 - \sqrt{7}\cos\left(\alpha - \frac{\pi}{3}\right)\right) = 0.80193\dots, \en 
where \eq{alpha}\alpha = \frac{1}{3}\arctan(3\sqrt{3}). \en
\end{thm} 
	
\begin{proof} It suffices to consider matrices of the form \eqref{n5r3}, with the parameters satisfying \eqref{bcdt}. The value of $\phi$ in \eqref{spec} is then inconsequential, and the multiple $e^{i\phi}$ can be ignored. 

We need to figure out when for $A$ given by \eqref{n5r3} and some unimodular $\omega$ the matrix $\re(\omega A)$ has multiple eigenvalues. Such an eigenvalue, if exists, has to also be an eigenvalue of the left upper 4-by-4 block of $\re(\omega A)$. 
Direct computations show that the eigenvalues of this block do not depend on $\omega$ and are as follows: \eq{lc} \lambda= \frac{1}{2}\sqrt{1\pm c},\ -\frac{1}{2}\sqrt{1\pm c}.\en  

The next step is to check when there exists a value of $\omega=e^{-i\theta}$ for which plugging $\lambda$ from \eqref{lc} into \eqref{n5r3kip-poly} yields zero. This happens if and only if the product of the extremal values 
$-\lambda^{5}  + \frac{3}{4}\lambda^{3} - \frac{1}{16} \left(c^{2} t^{2} + b^{2} + 1\right)\lambda - \frac{1}{16} bcst$ and 
$-\lambda^{5}  + \frac{3}{4}\lambda^{3} - \frac{1}{16} \left(c^{2} t^{2} + b^{2} + 1\right)\lambda + \frac{1}{16} bcst$ of \eqref{n5r3kip-poly} as a function of $\theta$, is non-positive. 

Equivalently, plugging in $\lambda$ as in \eqref{lc}:
\[ (1\pm c)\left( (1\pm c)^2-3(1\pm c)+1+b^2+c^2t^2\right)^2\leq 4 b^2c^2t^2s^2.\]

Substituting $b^2=1-c^2$, $s^2=1-t^2$ and solving for $t^2$ we see that the latter inequality holds only when $t=1/\sqrt{2\mp c}$, in which case it turns into the equality. The respective values of $\omega$ are $\pm 1$. 

Finally, for the such chosen $t$ the number $\pm\frac{1}{2}\sqrt{1+c}$ is a multiple eigenvalue of $\re A$ if and only if it is a root of 
\[ -5\lambda^4+\frac{9}{4}\lambda^2-\frac{1}{16}\left(2-c^2+\frac{c^2}{2-c}\right),\]
the derivative of \eqref{n5r3kip-poly}. This holds if and only if $c^3-2c^2-c+1=0$. As it happens, $c_+$ given by \eqref{cp} is the only such $c$ in $[0,1]$.

Similarly, $\pm\frac{1}{2}\sqrt{1-c}$ is a multiple eigenvalue of $\re A$ if and only if $c^3+2c^2-c-1=0$, and $c_-$ given by \eqref{cm} is the only such value in $[0,1]$. \end{proof} 	
Note that matrices \eqref{spec} both with $c_+$ and $c_-$ are non-generic. However, choosing $c_-$ yields a repeated non-extreme  eigenvalue of $\re(e^{-i\phi}A)$ and thus there is no flat portion on the boundary of the respective numerical range. On the other hand, for $c=c_+$ the eigenvalues of $\re(e^{-i\phi}A)$, counting the multiplicities, are approximately

$\{-0.75688, -0.36660,  -0.12348, 0.62348, 0.62348\}.$

An orthonormal basis of the eigenspace corresponding to the repeated eigenvalue can be chosen (also approximately) as 
\begin{equation*} \left\{ \begin{bmatrix}0.58217\\
0.02990 \\
0.72595 \\
0.21031\\
-0.29820\end{bmatrix}, \begin{bmatrix}0.00000\\
0.60104 \\
0.00000 \\
0.62348 \\
0.50000 \end{bmatrix}\right\}. 
\end{equation*}

The matrix of the compression of $\im(e^{-i\phi}A)$ in this basis is  

\eq{comp}\begin{bmatrix}
  0 & -i0.08077\dots\\
  i0.08077\dots & 0\\
\end{bmatrix}. \en

The preceding computations were carried out in Mathematica with 50 digits of precision. So, it seems that the flat portion actually materializes and has endpoints $e^{i\phi}(0.62349\dots\ \pm \ i0.08077 \dots)$. 

\begin{figure}[H]
\centering
    \includegraphics[scale=0.5]{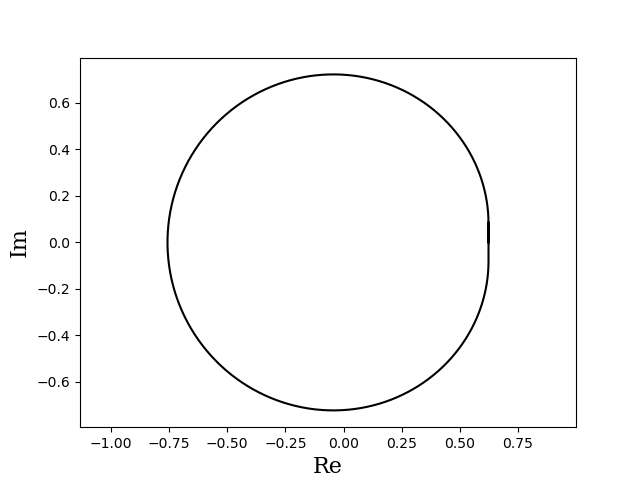}
    \caption{A plot of $W(A)$ for $c= c_+$ and $\phi = 0$.}
\end{figure}

\section{On higher rank numerical ranges}\label{s:hr} 
The rank-$k$ numerical range of $A\in\M_n$ is the set $\Lambda_k(A)$ defined as the set of $\lambda\in\C$ for which there exists an orthogonal projection $P\in\M_n$ of rank $k$ such that $PAP=\lambda P$ (see, e.g., \cite[Section 8.5]{GauWu}). It is clear from the definition that 
\[ W(A)=\Lambda_1(A)\supseteq\Lambda_2(A)\supseteq\ldots\supseteq \Lambda_n(A), \]
with the latter set non-empty only for $A=\lambda I$, in which case $\Lambda_n(A)=\{\lambda\}$. A much deeper property is that all the sets $\Lambda_k(A)$ are convex. More specifically, 
\eq{nrk} \Lambda_k(A)=\{ z\in\C\colon \lambda_{n-k+1}(\theta)\leq\re(e^{-i\theta}z)\leq\lambda_k(\theta), \theta\in(-\pi,\pi])\} \en
(see \cite[Theorem 5.11 in Chapter 8]{GauWu}). From \eqref{nil42poly} and the proof of Proposition~\ref{th:cir5} we therefore immediately obtain 
\begin{thm}\label{th:hnr}For $A$ unitarily similar to \eqref{nil42}, $\Lambda_2(A)$ is the circular disk centered at the origin of radius $r_2$ given by \eqref{r12} while $\Lambda_3(A)=\emptyset$.  
	
In turn, for $A$ unitarly similar to \eqref{n5r3} with condition \eqref{bcdx} satisfied, $\Lambda_2(A)$ is the circular disk centered at the origin of the radius $r_-$ given by \eqref{rpm}, $\Lambda_3(A)=\{0\}$, and $\Lambda_4(A)=\emptyset$. \end{thm}  

\providecommand{\bysame}{\leavevmode\hbox to3em{\hrulefill}\thinspace}
\providecommand{\MR}{\relax\ifhmode\unskip\space\fi MR }
\providecommand{\MRhref}[2]{%
	\href{http://www.ams.org/mathscinet-getitem?mr=#1}{#2}
}
\providecommand{\href}[2]{#2}

\end{document}